\documentstyle[12pt]{article}
\newcommand{\const}{\mathop{\rm const}\limits}

\newcommand{\Var}{\mathop{\rm Var}\limits}

\newcommand{\Law}{\mathop{\rm Law}\limits}

\begin{document}

\begin{center}

{\bf  Entropy and Grand Lebesgue Spaces approach \\

\vspace{3mm}

 for Prokhorov-Skorokhod continuity  of random  processes,   \\

\vspace{3mm}

with tail estimates.} \\

\vspace{4mm}

 {\sc Ostrovsky E., Sirota L.}\\

\vspace{3mm}

 \ Department of Mathematics and Statistics, Bar-Ilan University, \\
59200, Ramat Gan, Israel.\\
e-mail: \ eugostrovsky@list.ru \\

\vspace{3mm}

 \ Department of Mathematics and Statistics, Bar-Ilan University,\\
59200, Ramat Gan, Israel.\\
e-mail: \ sirota3@bezeqint.net \\

\vspace{3mm}

 {\sc Abstract.}

\end{center}

 \  We present in this paper a new sufficient condition for the  so-called Prokhorov-Skorokhod continuity
of random processes. Our conditions will  be formulated in the terms of metric entropy generated by three-dimensional
distribution of the considered random process (r.p.) in the parametric set, have a convenient and closed form, and generalize
some previous results. \par

 \ We study also the conditions for weak compactness of the sequence of random processes in this space and as a consequence
the Central Limit Theorem.\par

 \ Our consideration based on the theory of Prokhorov-Skorokhod spaces of random processes  and  on the
theory of  Banach spaces of random variables with exponential decreasing tails of distributions,
namely, on the theory of Grand Lebesgue Spaces (GLS) of random variables. \par

 \vspace{3mm}

  \par

\vspace{4mm}


\vspace{5mm}

\section{ Introduction. Previous works.} \par

\vspace{4mm}

 \hspace{3mm} Let $  f = f(t), \ t \in [0,1]  $ be real (or complex) valued  measurable function.
Recall that the  Prokhorov-Skorokhod module  $  \kappa[f](\delta)   $ for the function $  f(\cdot) $ at the point
$ \delta, \ \delta \in [0,1] $ is defined as follows:

$$
\kappa[f](\delta)  \stackrel{def}{=}
$$

$$
 \sup \{ \min|f(t) - f(t_1)|, \ |f(t_2) - f(t)|:
 0 \le t_1 \le t \le t_2 \le t_1 + \delta \le 1 \}.\eqno(1.1)
$$

 \ By definition, the function $ f: [0,1] \to R  $ belongs to the Prokhorov-Skorokhod  space  $  D[0,1]  $
iff

$$
\lim_{\delta \to  0+} \kappa[f](\delta) = 0 \eqno(1.2)
$$
and in addition

$$
\lim_{t \to 0+} f(t) = f(0), \hspace{5mm} \lim_{t \to 1-0} f(t) = f(1) \ - \eqno(1.3)
$$
unilateral continuity at both the boundary points $  t = 0 $ and $ t = 1.$ \par

\ Other name: "Cadlag" functions. These functions have in each point left and right limits; on the other words, without the
points of discontinuity of a second kinds, as in the classical book of I.I.Gikhman and A.V.Skorokhod, see \cite{Gikhman1}, chapters 4, 9.  \par
 \ We will take as ordinary that these functions are right continuous.

\vspace{3mm}

 \ There are (as a minimum) two version of the relation (1.3)  in the  case   when  $  f(t)  $ is a random process:
$  f(t) = \xi(t), \ t \in [0,1]:  $

$$
\lim_{ t \to 0+ } {\bf E} \arctan | \ \xi(t) - \xi(0) \ | = 0, \ \lim_{ t \to 1 - 0} {\bf E} \arctan | \ \xi(t) - \xi(1) \ | = 0, \eqno(1.3a)
$$
convergence in probability; \par

$$
{\bf P} ( \lim_{t \to 0+} \xi(t) = \xi(0)) = 1, \ {\bf P} ( \lim_{t \to 1 - 0} \xi(t) = \xi(1)) = 1, \  \eqno(1.3b)
$$
convergence almost everywhere. \par

\vspace{3mm}

  \ It is known that the function $  f: [0,1] \to R   $ belongs to the space  $  D[0,1]  $ iff it is right continuous,
has a left limits in each interior point, $  f(0+) = f(0), \ f(1-0) = f(1).  $ For instance, the trajectories of separable square
integrable martingales and empirical function of distribution  are elements of this space with probability one. \par
 \ The theory of this spaces is thoroughly outlined in the source articles and famous  books
 \cite{Billingsley1},  \cite{Billingsley2}, \cite{Gikhman1}, \cite{Kolmogorov1},  \cite{Prokhorov2},
\cite{Skorokhod1}; see also \cite{Bezandry1}, \cite{Bloznelis1}, \cite{Paulauskas1}, \cite{Sagitov1} and so one.\par
 In particular, this space  is complete metrizable relative appropriate distance  and  is herewith separable. \par

\vspace{3mm}

 \ Let a separable numerical valued random process $  \xi = \xi(t), \ t \in [0,1]  $ be given; for example,  is  given its
consistent family of finite-dimensional distributions.\par

\vspace{4mm}

 \ {\bf  Our target in this article is finding of simple closed sufficient conditions in the entropy and Grand Lebesgue Spaces (GLS)
terms for belonging of almost all paths of considered random process to the Prokhorov-Skorokhod  space:  }

$$
{\bf P} (  \xi(\cdot) \in D[0,1] ) = 1 \eqno(1.4)
$$
{\bf  and deriving the non-asymptotical estimates for tail distribution of the Prokhorov-Skorokhod module }

$$
T_{\Delta[\xi]}(u) = {\bf P} (\Delta[\xi]  > u). \eqno(1.5)
$$

 \ We  recall the needed definitions further, in the second and third sections. \par

 \ The exact exponential non-asymptotical estimates for tail  distribution of  an uniform norm  for the {\it random field} $  \xi(\cdot)  $

$$
 T_{||\xi||}(u) := {\bf P} ( \sup_{t \in [0,1]^d} |\xi(t)| > u), \ u \ge 1 \eqno(1.5a)
$$
are obtained, e.g. in the authors preprint \cite{Ostrovsky7}. \par

 \ Note that in contradistinction to the tail estimate in (1.5), the estimate in (1.5a)  is important not only at $ u \to \infty, $
but  also at $   u \to 0+. $ \par

\vspace{4mm}

 \ We will touch also briefly on the topic of weak compactness of the sequence $ \xi_n(\cdot)  $ of random processes in the space $  D[0,1]  $
with application to the Central Limit Theorem (CLT). \par

 \vspace{4mm}

 \ We will mention aside from the well-known applications  of the considered here problem
 in the theory of martingales and in the non-parametrical statistics also very interest applications, indeed: in physics \ \cite{Bloznelis2}, \cite{Daniels1}, and in the Monte-Carlo method \cite{Grigorjeva1} by computation of integrals from discontinuous integrand functions. \par

\vspace{4mm}

\section{Auxiliary notions and  facts.}

\vspace{4mm}

 \hspace{3mm} We present here for beginning some known facts from the theory of one-dimensional random variables
with exponential decreasing tails of distributions, see    \cite{Kozachenko1}, \cite{Ostrovsky1}, chapters 1,2. \par

 \ Especially  we menton the authors preprints \cite{Ostrovsky6}, \cite{Ostrovsky7};
 we offer in comparison with existing there results a more fine approach. \par

\vspace{3mm}
 \  Let $ (\Omega,F,{\bf P} ) $ be a probability space, $ \Omega = \{\omega\}. $ \par

\vspace{3mm}

 \ Let also $ \phi = \phi(\lambda), \lambda \in (-\lambda_0, \lambda_0), \ \lambda_0 =
\const \in (0, \infty] $ be certain even strong convex which takes positive values for positive arguments twice continuous
differentiable function, briefly: Young - Orlicz function,  such that
$$
 \phi(0) = \phi'(0) = 0, \ \phi^{''}(0) > 0, \ \lim_{\lambda \to \lambda_0} \phi(\lambda)/\lambda = \infty. \eqno(2.1)
$$

 \ For instance: $  \phi(\lambda) = 0.5 \lambda^2, \ \lambda_0 = \infty. $ \par

  \ We denote the set of all these  Young-Orlicz  function as $ \Phi; \ \Phi = \{ \phi(\cdot) \}. $ \par

 \ We say by definition that the {\it centered} random variable (r.v) $ \xi = \xi(\omega) $
belongs to the space $ B(\phi), $ if there exists some non-negative constant
$ \tau \ge 0 $ such that

$$
\forall \lambda \in (-\lambda_0, \lambda_0) \ \Rightarrow
\max_{\pm} {\bf E} \exp(\pm \lambda \xi) \le \exp[ \phi(\lambda \ \tau) ]. \eqno(2.2)
$$

 \ Obviously, this condition is quite equivalently to the well-known Kramer's condition

$$
\exists  \mu = \const > 0 \ \Rightarrow \max ( {\bf P}(\xi > x), \ {\bf P}(\xi < -x) ) \le \exp (-\mu x), \ x > 0. \eqno(2.2a)
$$

 \ The minimal non-negative value $ \tau $ satisfying (2.2) for all the values $  \lambda \in (-\lambda_0, \lambda_0), $
is named a $ B(\phi) \ $ norm of the variable $ \xi, $ write

$$
||\xi||B(\phi)  \stackrel{def}{=}
$$

 $$
 \inf \{ \tau, \ \tau > 0: \ \forall \lambda:  \ |\lambda| < \lambda_0 \ \Rightarrow
  \max_{\pm}{\bf E}\exp( \pm \lambda \xi) \le \exp(\phi(\lambda \ \tau)) \}. \eqno(2.3)
 $$

 \ These spaces are very convenient for the investigation of the r.v. having a
exponential decreasing tail of distribution, for instance, for investigation of the limit theorem,
the exponential bounds of distribution for sums of random variables,
non-asymptotical properties, problem of continuous and weak compactness of random fields,
study of Central Limit Theorem in the Banach space etc.\par

 \ The space $ B(\phi) $ with respect to the norm $ || \cdot ||B(\phi) $ and
ordinary algebraic operations is a rearrangement invariant Banach space which is isomorphic to the subspace
consisting on all the centered variables of Orlicz's space $ (\Omega,F,{\bf P}), N(\cdot) $
with $ N \ - $ function

$$
N(u) = \exp \phi^*(u) - 1, \hspace{4mm} \phi^*(u) \stackrel{def}{=} \sup_{\lambda} (\lambda u - \phi(\lambda)).
$$
 \ The transform $ \phi \to \phi^* $ is called Young-Fenchel transform. The proof of considered
assertion used the properties of saddle-point method and theorem of Fenchel-Moraux:
$$
\phi^{**} = \phi.
$$

 \ The next facts about the $ B(\phi) $ spaces are proved in \cite{Kozachenko1}, \cite{Ostrovsky1}, p. 19-40:

$$
 \xi \in B(\phi) \Leftrightarrow {\bf E } \xi = 0, \ {\bf and} \ \exists C = \const > 0,
$$

$$
U(\xi,x) \le \exp(-\phi^*(Cx)), x \ge 0, \eqno(2.4)
$$
where $ U(\xi,x)$ denotes in this section the {\it one-dimensional tail} of
distribution of the r.v. $ \xi: $

$$
U(\xi,x) = \max \left( {\bf P}(\xi > x), \ {\bf P}(\xi < - x) \right), \ x \ge 0,
$$
and this estimation is in general case asymptotically as $ x \to \infty  $ exact. \par

 \ Here and further $ C, C_j, C(i) $ will denote the non-essentially positive
finite "constructive" constants;  $ f^{-1}(\cdot) $ denotes the inverse function to the
function $ f $ on the left-side half-line $ (C, \infty). $ \par
  \ Let $  F =  \{  \xi(t) \}, \ t \in T, \ T  $ is an arbitrary set, be the family of somehow
dependent mean zero random variables. The function $  \phi(\cdot) $ may be "constructive" introduced by the formula

$$
\phi(\lambda) = \phi_F(\lambda) \stackrel{def}{=} \max_{\pm} \log \sup_{t \in T}
 {\bf E} \exp(  \pm \lambda \xi(t)), \eqno(2.5)
$$
 if obviously the family $  F  $ of the centered r.v. $ \{ \xi(t), \ t \in T \} $ satisfies the  so-called
{\it uniform } Kramer's condition:
$$
\exists \mu \in (0, \infty), \ \sup_{t \in T} U(\xi(t), \ x) \le \exp(-\mu \ x),
\ x \ge 0.
$$
 In this case, i.e. in the case the choice the function $ \phi(\cdot) $ by the
formula (2.5), we will call the function $ \phi(\lambda) = \phi_0(\lambda) $
a {\it natural } function, and correspondingly the function

$$
\lambda \to {\bf E} e^{\lambda \xi}
$$
is named often as a moment generating function for the r.v.  $ \xi, $ if of course there exists  in some non-trivial
neighborhood of origin. \par

\vspace{3mm}

 \ Further,  define the function $ \psi(p) = p/\phi^{-1}(p), \ p \ge 2. $ \par

\vspace{3mm}

 \ Let us introduce a new norm, the so-called "moment norm", or equally Grand Lebesgue Space (GLS) norm,
on the set of r.v. defined in our probability space by the following way: the
space $ G(\psi) $ consist, by definition, on all the centered (mean zero) r.v. with finite norm

$$
||\xi||G(\psi) \stackrel{def}{=} \sup_{p \ge 1} \left[ |\xi|_p/\psi(p) \right],  \eqno(2.6)
$$
here and in what follows as ordinary

$$
\ |\xi|_p :={\bf E}^{1/p} |\xi|^p = \left[  \int_{\Omega} |\xi(\omega)|^p \ {\bf P}(d \omega)  \right]^{1/p}.
$$

 \ It is proved that the spaces $ B(\phi) $ and $ G(\psi) $ coincides: $ B(\phi) =
G(\psi) $ (set equality) and both
the norms $ ||\cdot||B(\phi) $ and $ ||\cdot|| $ are linear equivalent: $ \exists C_1 =
C_1(\phi), C_2 = C_2(\phi) = \const \in (0,\infty), \ \forall \xi \in B(\phi) \ \Rightarrow $

$$
||\xi||G(\psi) \le C_1 \ ||\xi||B(\phi) \le C_2 \ ||\xi||G(\psi). \eqno(2.7)
$$

 \ In particular, let $ \eta $ be a numerical mean zero r.v. and let $  m = \const > 1. $ The following assertions are equivalent:

$$
{\bf A.} \ \exists C_1 \in (0, \infty) \ \Rightarrow \ U(\eta,x) \le \exp(-C_1 x^m), \ x \ge 0.\eqno(2.8a)
$$

$$
{\bf B.} \ \sup_{p \ge 1} \left[ \frac{|\eta|_p}{p^{1/m}} \right] < \infty. \eqno(2.8b)
$$

$$
{\bf C.} \  \exists C_2 \in (0, \infty) \ \Rightarrow  {\bf E} \exp(\lambda \eta) \le \exp \left( C_2 |\lambda|^{m/(m-1)} \right), \
|\lambda| \ge 1. \eqno(2.8c)
$$

  \ The definition (2.6) may be extended as follows. Recently, see \cite{Fiorenza1}, \cite{Fiorenza2}, \cite{Iwaniec1}, \cite{Kozachenko1},
\cite{Ostrovsky1}, chapters 1,2;  \cite{Ostrovsky7},  \cite{Rogoveer1}  appears the so-called Grand Lebesque Spaces (GLS) $  G(\psi) =
G(\psi; b) $ spaces consisting on all the measurable functions (random variables) $ \xi: \Omega \to R $ with finite norms

     $$
     ||\xi||G(\psi) \stackrel{def}{=} \sup_{p \in (1,b)} \left[ |\xi|_p /\psi(p) \right].
     $$

 \  Here $ \psi(\cdot) $ is some continuous positive on the open interval
    $ (1,b), \ b = \const \in (1, \infty]  \ $ function such that

     $$
     \inf_{p \in (1,b)} \psi(p) > 0.
     $$
  \  It is evident that $ G(\psi; b) $ is a rearrangement invariant space. \par

  \ Let  now $  l  $ be arbitrary number from the interval $ [1, \infty). $ We define a so-called {\it degenerate}
 GLS as follows. Put

 $$
 \psi_{(l)}(p) = \infty, \ p \ne l; \ \psi_{(l)}(l) = 1
 $$
and define formally $  \const/\infty = 0. $ Then the $ G\psi_{(l)}  $ norm of arbitrary r.v. $  \xi  $ coincides with the classical
$   L_l $ its norm:

$$
||\xi|| G\psi_{(l)}  = \sup_{p \ge 1} \left[ \frac{|\xi|_p}{\psi_{(l)}(p)} \right] = |\xi|_l,
$$
 if of course there exists. Thus, the classical Lebesgue - Riesz spaces $ L_p  $ are particular, more precisely, extremal case of
Grand Lebesgue Spaces.\par
   \ These spaces are used, for example, in the theory of probability,
 theory of PDE, functional analysis, theory of Fourier series, theory of martingales etc.\par

 \vspace{3mm}

 \ Let us consider more detail example. The inequality of a form

$$
|\xi|_p \le C_1 p^{1/m} \ \log^s p, \ p \ge 2, \ C_1 = \const \in (0, \infty) \eqno(2.9)
$$
where $  m = \const > 0, \ s = \const \in R $ is completely equivalent to the tail estimate

$$
U(\xi, x) \le \exp \left(  - C_2(m,C_1) \ x^m \ (\ln x)^{-ms} \right), \ x \ge e. \eqno(2.9a)
$$

 \ See, for example, \cite{Ostrovsky12}; \cite{Ostrovsky1}, chapter 1.8, theorem 1.8.1. \par

\vspace{4mm}

 \ The theory of multidimensional $  B(\phi) =  B(\vec{\phi}) $ spaces and correspondingly
 of multidimensional $  G(\psi) =  G(\vec{\psi}) $ ones  is represented in the recent articles
 \cite{Ostrovsky6},  \cite{Ostrovsky7}; it is quite analogous to the explained one. \par

 \vspace{3mm}

 \section{Main results.}

 \vspace{3mm}

 \hspace{4mm} Let us return now to the introduced random processes, say $  \xi(t), \ t \in [0,1]. $ Define for this r.p. and for the
fixed (non-random) triplet of numbers $ (r,s,t)  $  for which $  0 \le r \le s \le t \le 1 $

$$
\delta[\xi](r,s,t) \stackrel{def}{=} \min(|\xi(s) - \xi(r)|, |\xi(t) - \xi(s)|),\eqno(3.1)
$$

 \ We will denote the set of all such a triplets  by $  R:  $

$$
R = \{ (r,s,t): \ 0 \le r \le s \le t \le 1 \}
$$
and define for the fixed value $ s \in (0,1) $ the appropriate set

$$
 R(s) =  \{ r,t:  0 \le r \le s, \ s \le t \le 1 \ \}.
$$

 \ Define also

$$
\Delta[\xi] \stackrel{def}{=} \sup_{s \in (0,1)} \sup_{(r,t) \in R(s) } \delta[\xi](r,s,t). \eqno(3.2)
$$

\vspace{3mm}

 \ {\bf Proposition 3.0.} \ If  (in our notations)

$$
\sup_{s \in (0,1)} {\bf P} (\delta[\xi](r,s,t)  > u) \le [G(t) - G(r)]^{\alpha} \ u^{- 2 \beta}, \eqno(3.3)
$$
where $ \alpha = \const > 1, \ \beta = \const > 0, \ u = \const > 0, $ and $ G: [0,1] \to R $ is continuous increasing
deterministic function, then

$$
{\bf P} (\Delta[\xi]  > u) \le K(\alpha,\beta) \ u^{-2 \beta} \ [G(1) - G(0)]^{\alpha}, \ K(\alpha,\beta) < \infty \eqno(3.4)
$$
and correspondingly

$$
{\bf P}(\kappa[\xi](h) > u) \le 2 \ K(\alpha,\beta) \ u^{-2 \beta} \ [G(1) - G(0)]^{\alpha} \  (\omega[G](2 h) )^{\alpha - 1}, \eqno(3.5)
$$
where $ \omega[G](h)  $ is ordinary module of continuity of the function $ G(\cdot): $

$$
\omega[G](h)  \stackrel{def}{=} \sup \{ | \ G(t) - G(r) \ |: \ r,t \in [0,1], |r - t| \le h \}, \ h \in [0,1].
$$

 \ If in addition the r.p. $  \xi(t) $ satisfies the "boundary" conditions (1.3a), then almost all the trajectories of  $ \xi(t) $
belongs to the space $  D[0,1]. $

\vspace{3mm}

 \ This statement is proved, for example, in the preprint \cite{Sagitov1}, lemma 7.15; see also  \cite{Billingsley1}, chapters 2,3. \par

\vspace{3mm}

\ {\bf Remark 3.1.} We will repeatedly  apply the following obvious extension of proposition 3.0. Suppose the inequality (3.3) is true
for some {\it set} of the values $ (\alpha,\beta)  $ for which

$$
 D \subset \{ (\alpha, \beta): \ \alpha > 1, \ \beta > 0 \},
$$
with at the same function $  G(\cdot).  $ Then

$$
{\bf P} (\Delta[\xi]  > u) \le
 \inf_{(\alpha,\beta) \in D} \left\{ K(\alpha,\beta) \ u^{-2 \beta} \ [G(1) - G(0)]^{\alpha} \right\}, \eqno(3.4a)
$$

$$
{\bf P}(\kappa[\xi](h) > u) \le 2
\inf_{(\alpha,\beta) \in D} \left\{ \ K(\alpha,\beta) \ u^{-2 \beta} \ [G(1) - G(0)]^{\alpha} \  (\omega[G](2 h) )^{\alpha - 1} \right\}. \eqno(3.5a)
$$

 \ Of course, one can that the function $ G(\cdot) $ also dependent on the parameters $ (\alpha, \beta) $ from this set. \par

\vspace{5mm}

 \ {\bf Lemma 3.1.} The "constant"  $ K(\alpha,\beta)  $ in (3.4), (3.5) allows the following estimate

$$
K(\alpha,\beta) \le \frac{ \left( 1 - 2^{ (1 - \alpha)/(4 \beta)} \right)^{-2 \beta}}{2^{(\alpha - 1)/2} - 1}
\stackrel{def}{=} \overline{K}(\alpha, \beta).\eqno(3.6)
$$

 \vspace{3mm}

 \ {\bf Proof.} Serik Sagitov in  \cite{Sagitov1}  proved that

$$
\forall \theta \in \left( 2^{(1-\alpha)/(2 \beta)},1 \right) \ \Rightarrow  K(\alpha,\beta) \le
\frac{2^{(1 - \alpha)/(2 \beta)} \ \theta^{-2 \beta} \ ( 1 - \theta)^{-2 \beta}}{1 - 2^{1 - \alpha} \theta^{- 2 \beta} \ }, \eqno(3.7)
$$
therefore

$$
K(\alpha,\beta) \le \inf_{\theta \in \left( 2^{(1-\alpha)/(2 \beta)},1 \right)} \
\left[\frac{2^{(1 - \alpha)/(2 \beta)} \ \theta^{-2 \beta} \ ( 1 - \theta)^{-2 \beta}}{1 - 2^{1 - \alpha} \theta^{- 2 \beta} \ } \right]. \eqno(3.7a)
$$
 The required estimate (3.6) can be obtained after substituting

$$
\theta_0 = 2^{ (1 - \alpha)/(4 \beta)}
$$
 into the right-hand side of inequality (3.7). Note that this value $  \theta = \theta_0 $ is asymptotical as $ \alpha \to 1+0$
optimal. \par
 \ Note in addition that as $  \alpha \to 1 + 0 \ \Rightarrow   $

$$
 \overline{K}(\alpha, \beta) \sim 2^{4 \beta + 1} \ \beta^{2 \beta} \ (\ln 2)^{ - 2 \beta - 1 } \ (\alpha - 1)^{ - 2 \beta - 1 }.
$$

 \vspace{4mm}

 \ The  claim of this section is obtaining  the tails estimated for $ \Delta[\xi], \hspace{5mm}  \kappa[\xi](h) $ in the terms of
Grand Lebesgue Space norm for $  \delta(r,s,t) $  and consequently via generated by its metric entropy. \par

 \ To be more precise, we introduce the following semi-distance on the set $ [0,1] $ by means of some $  \psi \ - $ function:

$$
\rho_{\psi}(r,t) \stackrel{def}{=} \sup_{s \in (0,1)} ||\delta(r,s,t)||G\psi \eqno(3.8)
$$
and as a particular case

$$
\rho_{\psi_{(p)}}(r,t)  \stackrel{def}{=} \sup_{s \in (0,1)} |\delta(r,s,t)|_p.\eqno(3.9)
$$
 Herewith $ r < t $ and $  s \in [r,t]; $ we define in the case $  t < r \ d_p(t,r): = d_p( r,t) $  and
$  \rho_{\psi}(t,r): =  \rho_{\psi}(r,t),  $ so that $ d_p(t,r) = d_p( r,t)  $ and $ \rho_{\psi}(t,r) =  \rho_{\psi}(r,t). $
For reasons of continuity $ d_p(t,t): =  0 $ and $ \rho_{\psi}(t,t): = 0. $  \par

\vspace{4mm}

 \ We intend to deduce in this section the sufficient conditions for the Prokhorov-Skorokhod continuity  of the r.p. $ \xi(t) $
in the terms of metric entropy   of the set $  [0,1] $ generated by the semi-distance functions  $ \rho_{\psi}(t,r)  $ and
$ d_p(t,r), $ likewise the problem of natural continuity of the r.p., see e.g. \cite{Buldygin3}, \cite{Dudley1},
 \cite{Fernique1}, \cite{Fernique2}, \cite{Kozachenko1}, \cite{Ostrovsky1}, chapter 3, \cite{Pizier1}. \par

\vspace{4mm}

 \ Recall for reader convenience that the so-called {\it covering numbers} $  N(X,q,\epsilon)   $ of a (compact)  metric space $ (X,q)  $
relative the semi-distance function $  q = q(x_1, x_2), \ x_{1,2} \in X  $
are defined as a minimal amount  of closed balls

$$
 B(x_i,q, \epsilon) = B(x_i,\epsilon)  \stackrel{def}{=} \{ y, \ y \in X, q(x_i, y) \le \epsilon \}, \eqno(3.10)
$$
which cover all the set $   X: $

$$
 N(X,q,\epsilon) = \min\{N:  \exists x_j, \ j = 1,2, \ldots,N: \ \cup_{j=1}^N B(x_j, q, \epsilon) = X \}. \eqno(3.11)
$$

 \ The {\it metric entropy} $  H(X,q,\epsilon) $ of a (compact)  metric space $ (X,q)  $ is by definition  the natural logarithm
 of the covering number:

$$
  H(X,q,\epsilon) := \ln N(X,q,\epsilon).
$$

\vspace{3mm}

 \ {\bf Remark 3.2.} At the same definition may be used still in the case when the
function $  q = q(x_1, x_2)  $ is symmetrical, non negative, but  does not satisfy in general case the triangle inequality,  for
example, when

$$
q(t,r) = |t - r|^{\alpha}, \ t,r  \in [0,1], \ \alpha = \const > 1.
$$

\vspace{3mm}

 \ Denote also $ \Theta =  \{ \vec{\theta}  \}, \ \varepsilon  = \{ \vec{\epsilon} \},   $

$$
\vec{\theta} = \theta =
 \{ \theta(1), \theta(2), \ldots, \theta(k), \ldots \}: \ \theta(j) > 0, \ \sum_{k=1}^{\infty}\theta(k) = 1;
$$

$$
\vec{\epsilon} = \epsilon =
 \{\epsilon(1), \epsilon(2), \ldots, \epsilon(k), \ldots \}:  \ \epsilon(1) = 1, \ k \to \infty \ \Rightarrow
\epsilon(k) \downarrow 0.
$$

\vspace{3mm}

 \ The next fact is a slight and closed generalization of a statement 7.15 in the article \cite{Sagitov1}. \par

\vspace{3mm}

{\bf  Proposition 3.1. } Suppose  the considered separable random process $ \xi = \xi(t), \ t \in [0,1]  $ is such that
for some symmetrical and non negative numerical function $  q = q(t,r), \ t,r \in (0,1)  $
and for strictly increasing positive numerical function $  \lambda = \lambda(u), \ u > 0 $ for which

$$
\lim_{u \to \infty} \lambda(u) = \infty
$$
 there holds

$$
\sup_{s \in (0,1)} \sup_{(r,t) \in R(s)} {\bf P} (\delta[\xi](r,s,t) > u) \le \frac{q(r,t)}{\lambda(u)}.\eqno(3.12)
$$
 Then

$$
{\bf P} (\Delta[\xi] > 2u) \le Q(q(\cdot), \lambda(\cdot); \ u), \eqno(3.13a)
$$
where

$$
Q(q(\cdot), \lambda(\cdot); \ u) \stackrel{def}{=}
 \inf_{ \{\epsilon(k)\} \in \varepsilon } \inf_{ \{\theta(k) \} \in \Theta}
\sum_{k=1}^{\infty} N([0,1], q, \epsilon(k+1)) \cdot \frac{\epsilon(k)}{\lambda(u \cdot \theta(k))}, \eqno(3.13b)
$$
if of course the right-hand side of estimate (3.13b) tends to zero as $  u \to \infty. $  \par

  \ Denote also

 $$
 \sigma[q](h) = h^{-1} \ \sup_{ (r,t): |r - t| \le 2h} q(r,t).
 $$
  \  It is assumed that
$$
\lim_{h \to 0+}   \sigma[q](h) = 0.
$$

 \ We assert also under these conditions

$$
{\bf P} (  \kappa[\xi](h) > u) \le  Q(q(\cdot), \lambda(\cdot); \ u) \cdot \sigma[q](h).
$$

 \ If in addition the r.p. $  \xi(t) $ is continuous at the extremal points $  t = 0, \ t = 1 $  in the sense (1.3a),
then the random process $ \xi(\cdot) $ belongs to the Prokhorov-Skorokhod space $ D[0,1]  $ with probability one.\par

\vspace{3mm}

{\bf Example 3.1.} Suppose $ \lambda(u) = u^{2 \beta}, \ \beta = \const > 0. $  The conditions of proposition 3.1.
are satisfied if $ N([0,1], q,\epsilon) \le C \ \epsilon^{-\gamma}, \ \epsilon \in (0,1),  $  where $ \gamma = \const < 1.  $\par
 Namely, it is sufficient to choose  in (3.13) the values

$$
\epsilon(k) = s^{k - 1}, \ \theta(k) = (1-\theta) \ \theta^k,
$$
where

$$
0 < s, \theta < 1, \ s^{1 - \gamma} < \theta^{2 \beta},
$$
wherein

$$
{\bf P} (\Delta[\xi] > 2u) \le \frac{L(\beta,\gamma,\theta,C,s)}{u^{2 \beta}}, \ u > 0. \eqno(3.14)
$$

\vspace{3mm}

{\bf Example 3.2.} Suppose alike the example 3.1 $ \lambda(u) = u^{2 \beta}, \ \beta = \const > 0. $ We retain also the condition
 (3.12). But we suppose now

$$
 N([0,1], q, \epsilon) \le C \ \epsilon^{-1} \ |\ln \epsilon|^{-\gamma_1}, \ \epsilon \in (0,1/e),
$$
where $  \gamma_1 = \const > 1. $  If we choose in the estimate (3.13b)

$$
\theta(k) = C(\nu) \ k^{-\nu}, \ \nu = \const > \max(1, (\gamma_1 - 1)/(2 \beta) )
$$
and $  \epsilon(k) = \exp(-k + 1), $ then there holds as before

$$
{\bf P} (\Delta[\xi] > 2u) \le \frac{L_1(\beta,\gamma_1, C, \nu)}{u^{2 \beta}}, \ u > 0. \eqno(3.14b)
$$

 \ The last example does not be obtained from the proposition 7.15 in \cite{Sagitov1}. \par

\vspace{4mm}

\section{ Constructive building of distance function. \
Grand Lebesgue Spaces approach.}

\vspace{4mm}

 \ We discuss in this section the possibility of constructive building of the functions $ q(\cdot, \cdot), \ \lambda(\cdot)  $
for the condition  (3.12) in the terms of the source random process $ \xi = \xi(t) $ and following through the random process (field)
 $  \delta = \delta(r,s,t). $ \par

 \ We recall for beginning the analogous natural approach for the problem of continuity  relative  appropriate distance function
and finding of the tail estimation for the maximum distribution for the random process (field) $  \eta(t), \ t \in T, $
where $  T $  is arbitrary set, see \cite{Fernique2}, \cite{Kozachenko1}, \cite{Ostrovsky1}, chapters 3,4; \cite{Pizier1} etc.\par

 \ Introduce the natural function for the r.p. $  \eta(t) $

$$
\gamma(p) := \sup_{t \in T} |\eta(t)|_p,
$$
and suppose its finiteness at last for one value $ p $ greatest than one. Then $ \gamma(\cdot) $ is some $  \psi \ -  $ function and
one can to construct the semi-distance function $  d_{\gamma}(r,t)  $ also by natural way

$$
d_{\gamma}(r,t) := ||\eta(r) - \eta(t)||G\gamma.
$$
or equivalently

$$
\tilde{d}(r,t) := || \eta(r) - \eta(t)||B\phi
$$
with appropriate Young-Orlicz exponential type function $  \phi = \phi(\lambda).  $\par

 \ For instance, in $  \eta(t) $ is (separable) Gaussian random field, then  one can choose  $ \phi(\lambda) = 0.5 \ \lambda^2  $ and

$$
\tilde{d}(r,t) = \sqrt{ \Var (\eta(r) - \eta(t))}
$$
is the so-called Dudley's distance. \par

 \ The majority of results in described above problem were obtained  in the terms of metric entropy $ H(S,d_{\gamma}, \epsilon)  $
 of arbitrary subsets $  S \subset T $  generated by the distance $ d_{\gamma}(r,t). $\par
 \ This approach was introduced by  R.M.Dudley, see e.g. \cite{Dudley1}, and X.Fernique  \cite{Fernique1} - \cite{Fernique4};
it was applied in particular to the problem of Central Limit Theorem in the space of continuous functions. \par

 \ An another approach is closely related with the too modern notion "majorizing measures",  see \cite{Fernique1} - \cite{Fernique4},
\cite{Rogoveer1}, \cite{Talagrand1} - \cite{Talagrand4}. \par

 \ We return now to the formulated above problem  for (discontinuous, in general case)  random process $ \xi = \xi(t). $
  Let us introduce the natural function for the r.p. $ \delta = \delta(r,s,t) $

$$
\nu(p) := \sup_{s \in (0,1)} \sup_{(r,t) \in R(s)} |\ \delta(r,s,t) \ |_p, \eqno(4.1)
$$
{\it and suppose its finiteness for some value } $  p_0 > 2. $   \ It is not excluded herewith that
$  \nu(b) < \infty, $ where as before  $  b = \sup \{ p, \ \nu(p) < \infty \}, $ if of course  $  2 < b < \infty. $ \par

 \ The function $  \nu = \nu(p)  $ belongs to the set $  G\Psi = G\Psi(b) $ and herewith

$$
\sup_{s \in (0,1)} \sup_{(r,t) \in R(s)} || \ \delta(r,s,t) \ ||G\nu = 1. \eqno(4.2)
$$
 \ The pseudo - distance $  w = w(r,t)  $ may be defined by the formula

$$
w(r,t) := \sup_{s \in (0,1)} || \ \delta(r,s,t) \ ||G\nu, \ 0 \le r < t \le 1,\eqno(4.3)
$$
and we put by definition

$$
w(r,t) :=  w(t,r), \ r > t; \hspace{5mm}  w(r,r) = 0; \eqno(4.3a)
$$
so that

$$
\sup_{s \in (0,1)} | \ \delta(r,s,t) \ |_p \le w(r,t) \ \nu(p) \eqno(4.4)
$$
and by virtue of Tchebychev's inequality

$$
\sup_{s \in (0,1)}  {\bf P} (\delta(r,s,t) > u) \le \frac{\nu^p(p) \cdot w^p(r,t)}{u^p}, \ u > 0, \ p \in [2,b). \eqno(4.5)
$$

\ We intend now to use the statement of proposition (3.0) taking into account Remark 3.1.
 Namely, we suppose that the (deterministic!) function
$ \ w(r,t) \  $ is continuous and hence it allows \ \cite{Ostrovsky8} \ an  estimation of a form

$$
w(r,t) \le | \ G(t) - G(r)  \ |, \eqno(4.6)
$$
where $ G: [0,1] \to R $ is some continuous increasing deterministic function. We can and will suppose without loss of generality
$  G(0) = G(0+) = 0.   $ \par

 \ We substitute into inequality (3.3) the values $  \alpha = p, \  \beta = p/2 $  and recall that here $ p \ge 2: $

$$
\sup_{s \in (0,1)} {\bf P} (\delta[\xi](r,s,t)  > u) \le \nu^p(p) \ [G(t) - G(r)]^p \ u^{- p}, \eqno(4.7)
$$
then in accordance with proposition (3.0)

$$
{\bf P} (\Delta[\xi]  > u) \le K(p,p/2) \ \nu^p(p) \ [G(1)]^p \ u^{- p}, \ u > 0. \eqno(4.8)
$$
 \ It is easy to estimate

$$
K(p,p/2) \le 3^p, \ p \ge 2.
$$
 \ So, we obtained in fact the following statement. \par

 \vspace{3mm}

{\bf Proposition 4.1.} We deduce under formulated above in this section notations and conditions

$$
{\bf P} (\Delta[\xi]  > u) \le 3^p \ \nu^p(p) \ [G(1)]^p \ u^{- p}, \ u > 0. \eqno(4.9)
$$

  \ As a slight consequence:

$$
{\bf P} (\Delta[\xi]  > u) \le \inf_{p \in (1,b)} \left\{ 3^p \ \nu^p(p) \ [G(1)]^p \ u^{- p} \right\}, \ u > 0. \eqno(4.9a)
$$

\vspace{3mm}

 \ Evidently, the last estimations (4.9), \ (4.9a) are essentially non-improvable. \par

 \ Let us  estimate also in  the introduced terms and conditions the Prokhorov-Skorokhod module
$ \kappa[\xi] (h), \ h \in (0,1).  $ We apply again the inequality (3.5), in which  we substitute $ \alpha = p, \ \beta = p/2   $
and under at the same restriction $ p \ge 2:  $

$$
{\bf P}(\kappa[\xi](h) > u) \le 2 \ K(p,p/2) \ u^{-p} \ \nu^p(p)
\ [G(1) - G(0)]^p \  (\omega[G](2 h) )^{p - 1}. \eqno(4.10)
$$

\ We proved in fact the following estimate. \par

\vspace{3mm}

{\bf Proposition 4.2.} We deduce under formulated above in this section notations and conditions

$$
{\bf P}(\kappa[\xi](h) > u) \le 2 \inf_{p \in [2,b)}
\left[ \ 3^p \ u^{-p} \ \nu^p(p) \ [G(1) - G(0)]^p \  (\omega[G](2 h) )^{p - 1} \right]. \eqno(4.11)
$$

 \ Evidently, under these conditions

$$
\forall \epsilon > 0 \ \Rightarrow  \lim_{h \to 0+} {\bf P}(\kappa[\xi](h) > u) = 0,\eqno(4.12)
$$
so that if as before in addition the r.p. $  \xi(t) $ is (unilateral) continuous at the extremal points $  t = 0, \ t = 1 $  in the sense (1.3a),
then the random process $ \xi(\cdot) $ belongs to the Prokhorov-Skorokhod space $ D[0,1]  $ with probability one.\par

\vspace{3mm}

 \  We will use the last estimate further, in the seventh section. \par

\vspace{3mm}

\ {\bf Example 4.1.}  Suppose $  b = \infty $ and that

$$
\forall p \in [1, \infty) \ \Rightarrow
\nu(p) \le C_1 p^m, \ C_1, m = \const > 0;
$$
then we obtain the following {\it exponential decreasing} tail estimates

$$
{\bf P} (\Delta[\xi]  > u) \le \exp\left(-C_2(m, C_1, G(\cdot)) \ u^{1/m} \right), \ u \ge 1,
$$
and correspondingly for the values $ u \ge \left[ \omega[G](2h) \ |\ln \omega[G](2h)| \right]^{-m} $

$$
{\bf P}(\kappa[\xi](h) > u) \le 2 \ [ \omega[G](2 h) ]^{-1} \
\exp \left( -C_3(m, C_1, G(\cdot)) \ u^{1/m} \ \omega[G](2 h)  \right).
$$

\vspace{4mm}

\section{ Moment estimates for tail of minimum distribution.}

\vspace{4mm}

 \hspace{3mm} In order to  apply the results of the last section, we need to estimate the tail of distribution of minimum for the set
random variables. The exponential ones were received in  \cite{Ostrovsky6}, \cite{Ostrovsky7}.\par

  \ Let us consider at first a two-dimensional case $  d = 2. $ Indeed, let $  (\xi, \eta)  $ be a two-dimensional random vector.
Introduce a so-called binary absolute moment

$$
 \nu_{\xi,\eta}(p_1, p_2) = \nu(p_1, p_2) \stackrel{def}{=} {\bf E} \ |\xi|^{p_1} \ |\eta|^{p_2}, \ p_1, p_2   = \const > 0,
 \eqno(5.1)
$$
and correspondent pseudo-norm

$$
| \ (\xi,\eta) \ |_{p_1, p_2} :=   \left[ {\bf E} \ |\xi|^{p_1} \ |\eta|^{p_2} \right]^{1/(p_1 + p_2)} =
 \left[ \nu_{\xi,\eta}(p_1, p_2) \right]^{1/(p_1 + p_2)}. \eqno(5.2)
$$

 \  Denote also

$$
D(\xi,\eta) = \{ (p_1, p_2): \ \nu_{\xi,\eta}(p_1, p_2) < \infty  \}, \eqno(5.3)
$$

$$
T_{\xi, \eta} (u,v)  = {\bf P}( |\xi| > u, \ |\eta| > v ) = {\bf P}( (|\xi| > u) \cap \ (|\eta| > v) ), \ u,v > 0,  \eqno(5.4)
$$
be a two-dimensional tail function for the random vector $  (\xi, \eta). $  We deduce alike the proof of Tchebychev's inequality

$$
\nu_{\xi,\eta}(p_1, p_2) = \int_{\Omega}|\xi|^{p_1} \ |\eta|^{p_2} \ {\bf P}(d \omega) \ge
$$

$$
\int_{|\xi| > u, \ |\eta|> v} |\xi|^{p_1} \ |\eta|^{p_2} \ {\bf P}(d \omega) \ge
\int_{|\xi| > u, \ |\eta|> v} u^{p_1} \ v^{p_2} \ {\bf P}(d \omega) = u^{p_1} \ v^{p_2} \ T_{\xi, \eta} (u,v).
$$
 \ Therefore,

$$
T_{\xi, \eta} (u,v) \le \frac{\nu_{\xi, \eta} (p_1,p_2)}{u^{p_1} \ v^{p_2}}; \hspace{4mm} \ u,v > 0. \eqno(5.5)
$$

 \ As a slight  consequence: introduce the following set

$$
D = D(\xi,\eta) := \{ (p_1, p_2): \ \nu_{\xi, \eta} (p_1,p_2) < \infty  \};
$$
then

$$
T_{\xi, \eta} (u,v) \le \inf_{ (p_1,p_2) \in D(\xi, \eta)}
 \left[ \ \frac{\nu_{\xi, \eta} (p_1,p_2)}{u^{p_1} \ v^{p_2}} \ \right]; \hspace{4mm} \ u,v > 0. \eqno(5.5a)
$$

\vspace{4mm}

 \ We obtain as a consequence, choosing  $  v = u > 0: $

$$
{\bf P}(\min( |\xi|, \ |\eta|) > u)  \le \frac{ {\bf E} \ |\xi|^{p_1} \ |\eta|^{p_2} }{u^{p_1 + p_2}}. \eqno(5.6)
$$

 \ If we take in addition $  p_1 = p_2 = p > 0,  $ then

$$
{\bf P}(\min( |\xi|, \ |\eta|) > u) \le \frac{ {\bf E} \ |\xi \ \eta|^p}{u^{2p}}. \eqno(5.6a)
$$

\vspace{3mm}

 \ We can transform the estimate (5.5a) as follows. Let us extend the function $ \nu_{\xi,\eta}(p_1, p_2)  $ on the whole plane
of the values $ (p_1, p_2): $

$$
 \nu_{\xi,\eta}(p_1, p_2) := + \infty, \ (p_1, p_2) \notin D(\xi,\eta),
$$
then we have for the values $ u >1, \ v > 1  $

$$
T_{\xi, \eta} (u,v) \le \inf_{ (p_1,p_2) \in R^2}
 \left[ \ \frac{\nu_{\xi, \eta} (p_1,p_2)}{u^{p_1} \ v^{p_2}} \ \right] =
$$

$$
\inf_{ (p_1,p_2) \in R^2} \exp \left( - p_1 \ln u - p_2 \ln v + \ln \nu_{\xi, \eta} (p_1,p_2) \right) =
$$

$$
 \exp \left( -\sup_{(p_1,p_2) \in R^2}( p_1 \ln u + p_2 \ln v - \ln \nu_{\xi, \eta} (p_1,p_2) ) \right) =
$$

$$
 \exp \left( -  (\ln \nu_{\xi, \eta})^* (\ln u, \ln v) \right),
$$
where the notation $ f^*(\cdot,\cdot)  $ stands for the two - dimensional Young-Fenchel transform

$$
 f^*(\lambda_1,\lambda_2) \stackrel{def}{=} \sup_{x_1,x_2 \in R_2} [ x_1 \lambda_1 + x_2 \lambda_2 - f(x_1, x_2) ].
$$

\vspace{4mm}

 \ We list further some properties of the introduced pseudo-norm. \par

\vspace{3mm}

{\bf 0. } Note first of all that the expression for  $ |\ (\xi, \eta) \ |_{p_1, p_2}  $ does not represent in general case the really norm,
still for the values $  p_1 = p_2 = 1. $ For instance, the "unit sphere"

$$
S = \{ (\xi, \eta): \ |\ (\xi, \eta) \ |_{p_1, p_2} \le 1  \}
$$
is not convex set. \par

\vspace{3mm}

{\bf 1.} The functional $ (\xi, \eta) \to  |\ (\xi, \eta) \ |_{p_1, p_2}  $ is positive homogeneous of a degree 1:

$$
  |\ (\lambda \xi, \lambda \eta) \ |_{p_1, p_2} = |\lambda| \cdot |\ (\xi, \eta) \ |_{p_1, p_2}, \ \lambda = \const.
$$

\vspace{3mm}

{\bf 2.} Non-negativity:

$$
|\ (\xi, \eta) \ |_{p_1, p_2} \ge 0; \hspace{4mm}  |\ (\xi, \eta) \ |_{p_1, p_2} = 0 \ \Leftrightarrow  \xi \cdot \eta \stackrel{a.e.}{=} 0.
$$
 \ The relation  $  \xi \cdot \eta \stackrel{a.e.}{=} 0 $  is named often  as a {\it disjointness} of the random variables $  \xi  $  and $ \eta,$
especially for indicator functions. \par

 \vspace{3mm}

{\bf 3.}  A simple estimate.

$$
\nu_{\xi, \eta} (p_1,p_2)  = {\bf E} |\xi|^{p_1} \ |\eta|^{p_2} \le
\left[ {\bf E}|\xi|^{\alpha p_1 } \right]^{1/\alpha}  \cdot  \left[{\bf E} |\eta|^{\beta p_2} \right]^{1/\beta} =
|\xi|_{\alpha p_1}^{p_1} \cdot |\eta|_{\beta p_2}^{p_2}.
$$
We used the H\"older's inequality; here $ \alpha,\beta, p_1, p_2 = \const \ge 1, 1/\alpha + 1/\beta = 1.  $ Following

$$
|\ (\xi, \eta) \ |_{p_1, p_2} \le \inf \left\{ \ |\xi|_{\alpha p_1}^{p_1/(p_1 + p_2)} \cdot |\eta|_{\beta p_2}^{p_2/( p_1 + p_2) }: \
  \alpha,\beta \ge 1, 1/\alpha + 1/\beta = 1  \right\}.
$$

 \ For example,

$$
|\ (\xi, \eta) \ |_{p_1, p_2} \le  \ |\xi|_{2  p_1}^{p_1/(p_1 + p_2)} \cdot |\eta|_{2  p_2}^{p_2/( p_1 + p_2) }.
$$

 \ If in addition the r.v. $  \xi, \eta $ are independent, then $ {\bf E} |\xi|^{p_1} \ |\eta|^{p_2} =
 {\bf E} |\xi|^{p_1} \cdot {\bf E} |\eta|^{p_2} $  and hence

$$
|\ (\xi, \eta) \ |_{p_1, p_2} \le  \ |\xi|_{p_1}^{p_1/(p_1 + p_2)} \cdot |\eta|_{p_2}^{p_2/( p_1 + p_2) }.
$$

 \vspace{4mm}

{\bf 4.} Estimation of pseudo-norm for sum of random vectors. \par

$$
p_1, p_2 \ge 1 \ \Rightarrow
|\xi_1 + \xi_2, \eta_1 + \eta_2|_{p_1, p_2}  \le  2^{1 - 2/(p_1 + p_2)} \times
$$

$$
\inf_{\alpha, \beta> 0, \ 1/\alpha + 1/\beta = 1}
\left[\left(\xi_1|^{p_1}_{\alpha p_1} + |\xi_2|^{p_1}_{\alpha p_1} \right) \ \cdot \
 \left(\eta_1|^{p_2}_{\beta p_2} + |\eta_2|^{p_2}_{\beta p_1} \right) \right]^{ 1/(p_1 + p_2)}.
$$

\vspace{3mm}

 \ Proof.  Let $ \alpha, \ \beta $ be two fixed numbers such that $ \alpha, \beta> 0, \ 1/\alpha + 1/\beta = 1. $
 We observe applying the last estimate

$$
\nu_{\xi_1 + \xi_2, \eta_1 + \eta_2}(p_1,p_2) \le |\xi_1 + \xi_2|_{\alpha p_1}^{p_1} \cdot
 |\eta_1 + \eta_2|_{\beta p_2}^{p_2}.
$$

 \ The triangle inequality for the classical Lebesgue-Riesz  spaces $ L(p) $  together with an elementary inequality

$$
(a + b)^p \le  2^{p - 1} \ (a^p + b^p), \ a,b > 0, \ p \ge 1
$$
gives us  the required estimate. \par
 \ For example,

$$
|\xi_1 + \xi_2, \eta_1 + \eta_2|_{p_1, p_2}  \le  2^{1 - 2/(p_1 + p_2)} \times
$$

$$
\left[\left(\xi_1|^{p_1}_{2 \ p_1} + |\xi_2|^{p_1}_{2 \ p_1} \right) \ \cdot \
 \left(\eta_1|^{p_2}_{2 \ p_2} + |\eta_2|^{p_2}_{2 \ p_1} \right) \right]^{ 1/(p_1 + p_2)}.
$$

\vspace{6mm}

 \ It is not hard to generalize these propositions into the $ d \ - $ dimensional case.  Namely, let
$  \vec{\xi} = \{ \xi(1), \xi(2), \ldots, \xi(d) \}  $ be $ d \ - $ dimensional random vector and
$  \vec{p} = \{p(1), p(2), \ldots, p(d)  \}, \ \vec{u} = \{u(1), u(2), \ldots, u(d)  \}  $
be $ d \ - $ dimensional numerical  vectors with positive entries. Then the tail function for $ \vec{\xi}  $
may be easily estimated as follows

$$
T_{\vec{\xi}}(\vec{u}) \stackrel{def}{=} {\bf P} \left( \cap_{j=1}^d \{ |\xi(j)| > u(j) \} \right)  \le
\frac{ {\bf E}  \prod_{j=1}^d |\xi(j)|^{p(j)}}{\prod_{j=1}^d u(j)^{p(j)}}.
$$
 \ In particular,

$$
{\bf P} (\min_j |\xi(j)| > u) \le \frac{ {\bf E}  \prod_{j=1}^d |\xi(j)|^{p(j)}}{ u^{\sum_j p(j)}},
$$
and

$$
{\bf P} (\min_j |\xi(j)| > u) \le \frac{ {\bf E} \ \tau^{p}}{ u^{d \cdot  p}},
$$
where

$$
 \tau = \prod_{j=1}^d |\xi(j)|.
$$

 \ We note in continuation of this theme. Assume that the r.v. $  \tau $ belongs to some Grand Lebesgue Space
$  G\psi_b, \ b = \const > 1; $  for instance, one can choose $ \ \psi(\cdot)  \ $    as a natural function  for the r.v.
$ \tau: \ \psi(p) := |\tau|_p. $  We derive  for the values $  p  \in [1,b) $

$$
{\bf E} \ \tau^p \le \psi^p(p)
$$
and hence for the values $ u \ge 1 $

$$
{\bf P} (\min_j |\xi(j)| > u) \le \frac{ \psi^p(p)}{ u^{d \cdot  p}} =
\exp \left\{ - dp \ln u + p \ \ln \psi(p)   \right\},
$$

$$
{\bf P} (\min_j |\xi(j)| > u) \le  \inf_{p \in (1,\infty)} \exp \left\{ - dp \ln u + p \ \ln \psi(p)   \right\} =
$$

$$
\exp \left[ - \sup_{p} \left\{ pd \ln u - \psi_1(p) \right\} \right], \eqno(5.7a)
$$
where $ \psi_1(p) =  p \ \ln \psi(p), \ p \in [1,b) $  and $ \psi_1(p) =  + \infty $ when  $  p \notin [1,b).  $ \par
 \ The last term  may be expressed in turn  through the Young-Fenchel transform of the function $ \psi_1(p), $ as well:

$$
{\bf P} (\min_j |\xi(j)| > u) \le \exp \left[ - \psi^*_1 \{\ln (u^d) \} \right], \  u > 1. \eqno(5.7b)
$$

\vspace{3mm}

{\bf  Example 5.1: an application.}  Let again $  \xi(t), \ t \in [0,1] $ be a separable r.p.  We deduce

$$
{\bf P} (\delta(r,s,t) > u) \le \frac{ {\bf E} | \xi(r) - \xi(s) |^{p_1} \ |\xi(s)  - \xi(t)|^{p_2} }{ u^{p_1 + p_2}}, \hspace{4mm}
p_1, p_2, u > 0, \eqno(5.8)
$$
and as a particular case

$$
{\bf P} (\delta(r,s,t) > u) \le \frac{ {\bf E} | \xi(r) - \xi(s) |^p \ |\xi(s)  - \xi(t)|^p }{ u^{2p}}, \hspace{4mm} p,u > 0.
\eqno(5.9)
$$
and

$$
{\bf P} (\delta(r,s,t) > u) \le \inf_{p_1, p_2 > 0}
\left[\frac{ {\bf E} | \xi(r) - \xi(s) |^{p_1} \ |\xi(s)  - \xi(t)|^{p_2} }{ u^{p_1 + p_2}} \right], \hspace{4mm}
u > 0, \eqno(5.10)
$$

$$
{\bf P} (\delta(r,s,t) > u) \le \inf_{p > 0}\left[ \frac{ {\bf E} | \xi(r) - \xi(s) |^p \ |\xi(s)  - \xi(t)|^p }{ u^{2p}} \right],
\hspace{4mm} u > 0. \eqno(5.11)
$$

 \ The right-hand side of the  proposition (5.11)  may be estimated in turn  by means of Cauchy's inequality as follows.
 Denote by   $  d_p(s,t) $ the Pizier's (semi \ - ) distance

$$
d_p(s,t) \stackrel{def}{=} |\xi(s) - \xi(t)|_p = \left[ {\bf E} |\xi(s) - \xi(t)|^p   \right]^{1/p},  \eqno(5.12)
$$
then

$$
{\bf P}(\delta[\xi](r,s,t) > u) \le u^{-2p} \ \sqrt{ {\bf E} |\xi(r) - \xi(s)^{2p} \ {\bf E} |\xi(s) - \xi(t)|^{2p} } =
$$

$$
\frac{d_{2p}^p(r,s) \ d_{2p}^p(s,t)}{ u^{2p}}, \hspace{4mm} u > 0.\eqno(5.13)
$$

 Of course,

$$
{\bf P}(\delta[\xi](r,s,t)  > u) \le \inf_{p > 0} \left\{  \frac{d_{2p}^p(r,s) \ d_{2p}^p(s,t)}{ u^{2p}} \right\},
 \hspace{4mm} u > 0 \eqno(5.14)
$$
and following

$$
\sup_{s \in (0,1)}{\bf P}(\delta[\xi](r,s,t)  > u) \le
\sup_{s \in (0,1)}\inf_{p > 0} \left\{  \frac{d_{2p}^p(r,s) \ d_{2p}^p(s,t)}{ u^{2p}} \right\}, \hspace{4mm} u > 0 \eqno(5.15)
$$

\vspace{3mm}

 \ Let us impose the following condition on the r.p. $  \xi(\cdot): $

$$
d_p(r,s) \ d_p(s,t) \le Z(p) \cdot |V(t) - V(r)|^l, \ l = \const > 1/p  \eqno(5.16)
$$
for some function $ Z = Z(p) $ from the set $ [2,b), b = \const > 2 $ and for some continuous increasing bounded function
 $ V: [0,1] \to R $ for which $  V(0) = 0. $ We deduce  choosing

$$
\alpha = l p > 1, \ \beta = p, \ G(t) = Z^{1/l}(2p) \ V(t):
$$

$$
{\bf P} (\kappa[\xi](h) > u ) \le 2 \ K(l \ p, p) \ Z^{1/l}(2p) \ V^{l p}(1) \ u^{-2p} \ \left\{ \omega[Z](2h) \right\}^{lp - 1}, \eqno(5.17)
$$

$  u > 0, \ h \in [0,1], \ p \in [2, \min( b,1/l)) $ and hence $ {\bf P} (\kappa[\xi](h) > u ) \le $

$$
 \inf_{ p \in [2, \min( b,1/l))}
\left[ 2 \ K(l \ p, p) \ Z^{1/l}(2p) \ V^{l p}(1) \ u^{-2p} \ \left\{ \omega[Z](2h) \right\}^{lp - 1} \right]. \eqno(5.18)
$$

 \ Note that under formulated above assumptions

$$
\forall \epsilon > 0 \ \Rightarrow \lim_{h \to 0+} {\bf P} (\kappa[\xi](h) > \epsilon ) = 0. \eqno(5.19)
$$

 \vspace{3mm}

 \section{ About boundary restrictions.}

 \vspace{3mm}

 \ Let's turn our attention to the condition (1.3) for the random process $  \xi(t), \ t \in  [0,1].$ Question:
under what conditions (sufficient or necessary or sufficient and necessary conditions) on the distribution on the $ \xi(\cdot) $

$$
{\bf P} (\lim_{t \to 0+} (\xi(t) - \xi(0)) = 0 ) = 1 \eqno(6.1)
$$
or analogously

$$
{\bf P} (\lim_{t \to 1-0} (\xi(t) - \xi(1)) = 0 ) = 1. \eqno(6.1a)
$$

 \ Denote

$$
Z_0(\beta) = {\bf E} \arctan \sup_{t \in [0,\beta]} |\xi(t) - \xi(0)|,
$$

$$
Z_1(\beta) = {\bf E} \arctan \sup_{t \in [1 -\beta,1]} |\xi(t) - \xi(1)|, \ \beta = \const \in (0,1/2).
$$

\vspace{3mm}

 \ {\bf Proposition 6.1.} \par
  \ {\bf A.} The condition

$$
\lim_{\beta \to 0+} Z_0(\beta) = 0 \eqno(6.2)
$$
 is necessary and sufficient for the equality (6.1). \par

  \ {\bf B.} The condition

$$
\lim_{\beta \to 0+} Z_1(\beta) = 0 \eqno(6.2a)
$$
 is necessary and sufficient for the equality (6.1a). \par

 \vspace{3mm}

\ {\bf The proof } is quite analogously to one for a main result of the author's preprint
\cite{Ostrovsky10}  and may be omitted. \par

 \ But it is worth to note that if

$$
\lim_{h \to 0+} \kappa[\xi](h)  = 0
$$
 with probability one and

$$
\lim_{t \to 0+} (\xi(t) - \xi(0)) = 0,  \eqno(6.3)
$$

$$
\lim_{t \to 1-0} (\xi(t) - \xi(1)) = 0  \eqno(6.3a)
$$
in the sense of convergence in probability, or equally (here) in the sense of convergence in distribution, then
$  {\bf P} (  \xi(\cdot) \in D[0,1]) = 1,  $ see \cite{Sagitov1}. \par

 \vspace{3mm}

 \section{Conditions  for weak compactness. CLT in this space.}

 \vspace{3mm}

 \hspace{3mm} Let $  X(t); X_n(t), n = 1,2,\ldots, \ t \in [0,1]  $ be a sequence of separable random processes.
We will study in this section the problem of finding sufficient conditions for weak (in distribution) convergence in the Prokhorov-Skorokhod space

$$
\Law (X_n(\cdot)) \stackrel{D[0,1]}{ \to } \Law (X(\cdot)).  \eqno(7.1)
$$

  We will suppose in this section that all the finite-dimensional distributions of r.p. $  X_n(t) $ converges to ones for the
r.p. $  X(t). $ Assume also that the limit process $  X(t)  $ satisfies the boundary conditions (1.3a). \par

 \ We will study in this section the problem of finding sufficient conditions for weak (in distribution) convergence in the Prokhorov-Skorokhod space

$$
\Law (X_n(\cdot)) \stackrel{D[0,1]}{ \to } \Law (X(\cdot)).  \eqno(7.1)
$$

 \ Of course, we are forced to admit that $ X_n(\cdot), \  X(\cdot) $ are elements  of the space $ D[0,1] $ with probability one. \par

 \ We need to introduce some new notations.  Define the following  uniform natural $ \ G\Psi \ $ function

$$
\zeta(p) = \zeta [\{X_n\} ](p) \stackrel{def}{=} \sup_n \sup_{s \in (0,1)} \sup_{(r,t) \in R(s)} |\delta[X_n(\cdot)], (r,s,t) |_p, \eqno(7.2)
$$
and suppose its finiteness for the at last one value $ p = p_0 = \const > 2; $ denote as before
$ b = \sup \{ p: \ \zeta(p) < \infty \}; $ then $  b = \const \in (2, \infty].  $ \par
 \ There exists a continuous increasing function $ Q: [0,1] \to R  $ for which $ Q(0) = Q(0+) = 0  $ and

$$
\sup_n \sup_{s \in (0,1)} | \ \delta[X_n(\cdot)](r,s,t) \ |_p \le |Q(t) - Q(r)| \cdot \zeta(p). \eqno(7.3)
$$

\vspace{4mm}

{\bf Proposition 7.1.} We deduce on the basis of proposition 4.1
under formulated above in this section notations and conditions

$$
\sup_n {\bf P} (\Delta[X_n]  > u) \le \inf_{p \in [2,b)} \left\{ 3^p \ \zeta^p(p) \ [Q(1)]^p \ u^{- p} \ \right\}, \ u > 0; \eqno(7.4)
$$

$$
 \sup_n {\bf P}(\kappa[X_n](h) > u) \le 2 \inf_{p \in [2,b)}
\left\{  \ 3^p \ u^{-p} \ \zeta^p(p) \ [Q(1)]^p \  (\omega[Q](2 h) )^{p - 1} \right\}. \eqno(7.5)
$$

 \ Furthermore, if in addition all the considered r.p. $  X_n(t), $ including the
limiting random process $ X(t),  $ are continuous at the extremal points $ t = 0, \ t = 1  $
in the sense (1.3a), then the sequence of r.p. $ X_n(\cdot) $ converges at the r.p. $  X(\cdot) $  weakly in distribution in
the space $  D[0,1].  $ \par

\vspace{3mm}

\ {\bf Proof.} The estimates (7.4), (7.5) follows immediately from the proposition (4.1).  Both these inequalities together with
the convergence of r.v.  $  X_n(0) - X(0) \to 0, \ X_n(1) - X(1) \to 0 $ guarantee us the weak compactness of the distributions
$  X_n(\cdot) $  in the Prokhorov-Skorokhod space $  D[0,1]. $ \par
 \ Finally, the convergence of all the finite-dimensional distributions of r.p. $  X_n(t) $ to the ones for $  X(t) $ gives us
what is required, see e.g. \cite{Sagitov1}. \par

\vspace{3mm}

 \ Let us now turn as a capacity of the particular case to the study of the Central Limit Theorem in this space.\par

 \vspace{3mm}

 \ We recall here the classical definition of the CLT in  Prokhorov-Skorokhod (or more generally in arbitrary linear separable
topological) space.  Let $  \xi(t) = \xi_1(t) $ be centered (mean zero) separable random process with values in this space having finite
variance in weak sense. Let $ \ \xi_i = \xi_j(t), \ j = 1,2,\ldots \  $ be independent copies of $  \xi(t). $ Denote

$$
S_n(t) = n^{-1/2} \sum_{j=1}^n  \xi_j(t),
$$
and  let $ S_{\infty}(t) = S(t)  $ be centered separable Gaussian process with at the same covariation as $ \xi(t). $ \par

 \ It will be presumed that all the random processes $  \xi_j(t), \ S_{\infty}(t)  $ are defined at the same  sufficiently rich
probability space. \par

 \ By definition, the r.p. $ \xi(\cdot),  $ or equally the sequence of r.p. $  \{  \xi_j(\cdot) \}  $ satisfies the Central Limit Theorem
(CLT) in the space $ D[0,1], $  if all the considered r.p. $  \{  \xi_j(\cdot) \}, \ S_{\infty}(\cdot)  $  belong to this space almost surely
and if the sequence of the distributions $  S_n(\cdot)  $ converges weakly as $  n \to \infty $ to the distribution of the $  S_{\infty}(\cdot): $
for arbitrary bounded continuous functional $  F: D[0,1] \to R  $

$$
\lim_{n \to \infty} {\bf E} F(S_n) = {\bf E} F(S_{\infty}).
$$
 In particular,

$$
\lim_{n \to \infty} {\bf P} (  ||S_n|| > u) = {\bf P} ( ||S_{\infty}|| > u), \hspace{4mm} u > 0.
$$

 \ The latter circumstance is the basis not only in the non-parametrical statistics, but also in the Physic \cite{Daniels1}, and
in the Monte-Carlo method for computation of multiple integrals from the discontinuous functions, see  \cite{Frolov1},
\cite{Grigorjeva1}. \par

\vspace{4mm}

 \ Evidently, the finite - dimensional distributions of the r.p. $  S_n(\cdot)  $ converges as $  n \to \infty $ to ones for
the r.p. $ S_{\infty}(t); $  it remains only to ground the weak compactness of the correspondent distributions in the space $  D[0,1]. $ \par

\vspace{3mm}

 \ We retain the definitions and result of the example 5.1., especially the estimates  (5.12) and (5.13); recall only that
$ {\bf E} \xi(t) = 0 $ and $  p \ge 2.  $ \par

 \ We will use the famous Rosenthal's inequality, see  \cite{Rosenthal1}, \cite{Johnson1},  \cite{Ibragimov1} etc. in the following form.
Let  $ \tau, \ \tau_i = 1,2,\ldots $  be a sequence of i., i.d. mean zero  random variables  with finite $ p^{th} $ absolute moment. The
following estimate holds true:

$$
\sup_{n \ge 1} \left| n^{-1/2} \sum_{i=1}^n \tau_i  \right|_p \le K_R(p) \cdot |\tau|_p, \ p \ge 2,  \eqno(7.6)
$$
where $  K_R(p) $ is so-called  Rosenthal's "constant", more precisely,  function on $  p. $ It is known, see \cite{Ostrovsky11},
that

$$
K_R(p) \le C_A \cdot \frac{p}{\ln p}, \ p \ge 2, \eqno(7.6a)
$$
where $  C_A $ is an absolute constant, with the following its value $  C_A \approx 0.65349368 < 0.6535.  $ \par

\vspace{3mm}

 \ We will apply the inequality (7.6) to the sequence of differences $  \xi_j(t) - \xi_j(r)  $ with correspondent norms

$$
d_p[\xi](t,r) := |\xi(t) - \xi(r)|_p,
$$
and suppose as before its finiteness for the values $ \ p \in [2,b), \ $ where $  b = \const \in (2, \infty].  $

 \ We conclude

$$
\sup_n d_p[S_n](t,r) \le K_R(p) \cdot d_p[\xi](t,r). \eqno(7.7)
$$

 \ Define analogously the following  natural $ \ G\Psi \ $ function

$$
y(p) = y[\xi](p) \stackrel{def}{=} \sup_{s \in (0,1)} \sup_{(r,t) \in R(s)} |\delta[\xi(\cdot)], (r,s,t) |_p, \ p \in [2,b), \eqno(7.8)
$$
then

$$
Y(p) \stackrel{def}{=} \sup_n y[S_n](p) \le  K_R(p) \cdot y(p). \eqno(7.9)
$$

 \ There exists a continuous increasing function $ B: [0,1] \to R  $ for which $ B(0) = B(0+) = 0  $ and

$$
 \sup_{s \in (0,1)} | \ \delta[\xi(\cdot)](r,s,t) \ |_p \le |B(t) - B(r)| \cdot y(p), \eqno(7.10)
$$
 following

$$
\sup_n \sup_{s \in (0,1)} | \ \delta[S_n(\cdot)](r,s,t) \ |_p \le |B(t) - B(r)| \cdot K_R(p) \cdot y(p). \eqno(7.11)
$$

 \ It remains to use the  proposition 7.1. \par

\vspace{4mm}

{\bf Proposition 7.2.} We deduce on the basis of proposition 7.1
under formulated above in this section notations and conditions

$$
\sup_n {\bf P} (\Delta[S_n]  > u) \le \inf_{p \in [2,b)} \left\{ 3^p \ K_R^p(p) \ y^p(p) \ [B(1)]^p \ u^{- p} \ \right\}, \ u > 0; \eqno(7.12)
$$

$$
 \sup_n {\bf P}(\kappa[S_n](h) > u) \le 2 \inf_{p \in [2,b)}
\left\{  \ 3^p \ u^{-p} \  K_R^p(p) \ y^p(p) \ [B(1)]^p \  (\omega[B](2 h) )^{p - 1} \right\}. \eqno(7.13)
$$

\vspace{4mm}

 \ As a consequence, the r.p. $ \xi(\cdot),  $ or equally the sequence of r.p. $  \{  \xi_j(\cdot) \}  $ satisfies the Central Limit Theorem
in the space $ D[0,1]. $ \par

\vspace{4mm}

 \ {\bf Example 7.1.} Assume that

$$
y(p) \le C_1 \ p^{1/m} \ \ln^s p, \ p \ge 2, \ m = \const > 0, \ s = \const; \eqno(7.14)
$$
then we have the following non-asymptotical tail estimates: $ \sup_n {\bf P} (\Delta[S_n]  > u) \le $

$$
 \exp \left\{-C_2(C_1, m,s) \ u^{m/(m+1)} \ | \ \ln u \ |^{m(s-1)/(m + 1)} \right\}, \ u \ge e;  \eqno(7.15)
$$

$$
 \sup_n {\bf P}(\kappa[S_n](h) > u) \le 2 \  (\omega[B](2 h) )^{- 1} \times
$$

$$
 \exp \left\{ - C_3(C_1,m,s) \ [ \ u/ (\omega[B](2 h)) \ ]^{ m/(m + 1)} \  [ \ \ln ( u/(\omega[B](2 h))) \ ]^{ m(s-1)/(m+1)  } \right\},
$$
when

$$
u >  e \cdot \omega[B](2 h) \cdot |\ \ln \omega[B](2 h) \ |^{1 + 1/m}. \eqno(7.16)
$$

 \ In turn, the condition (7.14) may  be expressed in the terms of tail behavior for the r.p.  $  \delta(r,s,t), $ see (2.9) - (2.9a). \par

 \vspace{3mm}

 \section{Concluding remarks.}

 \vspace{3mm}

 \hspace{4mm} {\bf A.} It is interest by our opinion to generalize obtained results into the multidimensional case, i.e. into the space
$  D[0,1]^d. $ The non-asymptotical estimated for tail of uniform norm distribution for discontinuous random processes
are obtained in \cite{Ostrovsky7}. \par

\vspace{3mm}

{\bf B.} It is interest also a generalization  on the case when the sequence of r.p. $ \{  \xi_j(t) \}  $ forms on the index
$ \ j  \ $ a sequence of martingale differences relative appropriate filtration. \par

\vspace{3mm}

{\bf C.}  Perhaps, the applying of the more modern technic, indeed the so-called  majorizing measures, see
 \cite{Fernique1}, \cite{Talagrand1} - \cite{Talagrand4} one can give more exact estimated. \par

\vspace{3mm}

{\bf D.} More interest new examples of CLT in the Prokhorov-Skorokhod space with applications may be found in the articles
\cite{Bezandry1}, \cite{Bloznelis1} - \cite{Bloznelis3}, \cite{Grigorjeva1}, \cite{Paulauskas1} etc. \par

\end{document}